\documentclass[11pt]{amsart}
\usepackage{amssymb,amsthm,verbatim, latexsym}
\newtheorem{theorem}{Theorem}[section]  
\newtheorem{lemma}[theorem]{Lemma}       
 \newtheorem{corollary}[theorem]{Corollary}     
\newtheorem{proposition}[theorem]{Proposition} 
\newtheorem{remark}[theorem]{Remark}           
\newtheorem{definition}[theorem]{Definition}   
\numberwithin{equation}{section}
\newcommand{\nc}{\newcommand}
\nc{\R}{\ensuremath{\mathbb{R}}}
\nc{\os}{\ensuremath{\mathfrak {S}}}
\nc{\cl}{\ensuremath{\operatorname{cl}}}
\nc{\C}{\ensuremath{\mathbb{C}}}
\nc{\K}{\ensuremath{\mathbf{K}}}
\nc{\N}{\ensuremath{\mathbb{N}}}
\nc{\sgn}{\operatorname{sgn}}
\nc{\id}{\operatorname{id}}
\nc{\A}{\ensuremath{\mathcal {A}}}
\nc{\E}{\ensuremath{\mathcal {E}}}
\nc{\OS}{\operatorname{OS}}
\newcommand{\I}{\ensuremath{\Im}}
\begin{document}
\title[An Orlik-Solomon type algebra for matroids]
{An Orlik-Solomon type algebra for matroids with a fixed linear class of circuits}
\thanks{2000 \emph{Mathematics Subject Classification}: 
05B35, 05Cxx, 14F40, 32S22.\\
 \emph{Keywords and phrases}: arrangement of
 hyperplanes, 
 matroid, 
Orlik-Solomon algebra, biased graphs.
}
  \author[Raul Cordovil and David Forge]{Raul Cordovil and David Forge}
\address{
  Raul Cordovil\newline
Departamento de Matem\'atica,\newline 
Instituto Superior T\' ecnico \newline
 Av.~Rovisco Pais
 - 1049-001 Lisboa  - Portugal}
\email{cordovil@math.ist.utl.pt}
\urladdr{http://www.math.ist.utl.pt/~rcordov}
\address{
David Forge\newline
Laboratoire de Recherche en Informatique UMR 8623\newline 
Batiment 490
Universit\'{e} Paris Sud\newline 
91405 Orsay Cedex - France}
\email{forge@lri.fr}
\urladdr{http://www.lri.fr/~forge}
\thanks{The  authors' research  was
supported in part by the project  ``Algebraic Methods in  Graph Theory" 
approved by program Pessoa 2005.  The first  author's research  was also
supported by FCT (Portugal) through 
program POCTI and while visiting the Laboratoire de Recherche en Informatique
of the Universit\'e Paris-Sud.}

%%%%%%%%%%%%%%%%%
\begin{abstract} A family $\mathcal{C}_L$ of circuits of a matroid $M$ is 
a linear class if, given a modular pair of circuits in $\mathcal{C}_L$, 
any circuit contained in the union of the pair is also in $\mathcal{C}_L$.
The pair $(M, \mathcal{C}_L)$ can be seen as a matroidal generalization 
of a biased graph. 
We introduce and study an Orlik-Solomon type algebra
determined by  $(M, \mathcal{C}_L)$. If $\mathcal{C}_L$ is the set of all 
circuits of $M$ this algebra is  the Orlik-Solomon algebra of $M$.
\end{abstract}
\maketitle
\today
\section{Introduction}
%%%%%%%%%%%
Let $\A_{\C}=\{H_{1},\ldots, H_{n}\}$ be a central and essential arrangement of 
hyperplanes in $\C^d$ (i.e, such that  $\bigcap_{H_i\in \A_{\C}} H_i= \{0\}$). The manifold $\mathfrak{M}=\C^d\setminus\bigcup_{H\in
 \mathcal{A}_{\C}} H$  plays  an important
role in the 
 Aomoto-Gelfand multivariable theory 
 of  hypergeometric functions (see \cite{OT2} for a recent 
introduction from the point of  view of
arrangement theory). There is a  rank $d$
matroid $M:=M(\mathcal{A}_{\C})$ on the ground set 
$[n]$  canonically determined by  
$\mathcal{A}_{\C}$: a subset $D\subseteq [n]$\, 
is  a \emph{dependent set} of  
$M$ if and only if there are  scalars $\zeta_{i}\in \mathbb 
{C},\, i\in D$, not all nulls,
such that $\sum_{i\in D}\zeta_{i}\theta_{H_{i}}=0$, where $\theta_{{H}_i}\in (\C^d)^{*}$
 denotes a linear form such that $\mathrm{Ker}(\theta_{{H}_i})={H_i}$.
 \par
Let $M$ be a matroid and $M^\star$ be its dual. In the following, 
we suppose that the ground set of $M$ is $[n]:=\{1,2,\dotsc,n\}$ and 
its  rank function is denoted by $r_M$. The subscript $M$ in $r_M$ will often be omitted.
Let $\mathcal{C}=\mathcal{C}(M)$ be the family of circuits of $M$. 
Let $\mathbf{K}$ be a field and $E=\{e_1,\ldots,e_n\} $ 
be a finite set of order $n$.  Let 
$\bigoplus_{e\in E} \mathbf {K}e$ be the vector space over $\mathbf {K}$ 
of basis $E$ 
and $\mathcal{E}$ be the graded exterior algebra $\bigwedge 
\big(\bigoplus_{e\in E} \mathbf {K}e\big)$, i.e.,
$$\mathcal{E}:=\sum_{i=0} \E_i=\mathcal{E}_0(= \mathbf {K}) 
\oplus  \mathcal{E}_1(=
\bigoplus_{e\in E} \mathbf {K}e)  \oplus \dotsm 
\oplus\mathcal{E}_i(= \bigwedge^i(\bigoplus_{e\in E}
\mathbf {K}e))\oplus \dotsm.$$ 
For every linearly ordered subset $X=\{{i_1},\dotsc,  {i_m}\}
\subseteq [n],\, i_1<\dotsm < {i_m},$\, let $e_X$ be the monomial
$e_X:=e_{i_1}\wedge e_{i_2}\wedge \dotsm\wedge e_{i_m}$. By definition 
  set  $e_\emptyset=1\in \mathbf{K}$. 
Consider the map $\partial: \mathcal{E} \to  \mathcal{E}$,
 extended by linearity from  the ``differentials'', 
$\partial e_i =1$ for every $e_i\in E$, $\partial e_\emptyset=0$ and
$$\partial e_{_X}=\partial(e_{i_1}\wedge \dotsm \wedge
e_{i_m})=\sum(-1)^{j}e_{i_1}\wedge\dotsm \wedge
{e}_{i_{j-1}}\wedge
{e}_{i_{j+1}}\wedge\dotsm
\wedge e_{i_m}.$$ 
The (graded) {\em Orlik-Solomon
$\mathbf{K}$-algebra}
$\OS(M)$   of the matroid $M$
is the quotient   
$\mathcal{E} 
/\Im$ where $\Im$ denotes the (homogeneous)  two-sided ideal of 
$\mathcal{E}$ generated by
the set 
$$\big\{\partial e_C:\, C\in \mathcal{C}(M), |C|>1\big\}
\cup \big\{e_C: C\in \mathcal{C}(M), |C|=1\big\}$$
 or equivalently by the set
$$\big\{\partial e_C:\, C\in \mathcal{C}(M), |C|>1\big\}
\cup \big\{e_C: C\in \mathcal{C}(M)\big\}.$$
The de Rham cohomology algebra $H^{\bullet}\big(\mathfrak{M}
(\mathcal{A}_{\mathbb{C}}); \mathbf{K}\big)$  is
shown to be isomorphic to the  Orlik-Solomon $\mathbf{K}$-algebra of the 
matroid 
$M(\A_{\mathbb{C}}),$ see \cite{OS1, OT}. 
We refer  to  \cite{Falk} for a
recent  
discussion on the role of matroid theory in the study of 
Orlik-Solomon algebras. \par
%%%%%%%%%%%%%%%%%%%%%%%%%%%%%%%%%%%%%%%%%%%%%%%%%
\section{Linear class of circuits}
Given a family $\mathcal{C}$ of circuits of a matroid $M$ set 
$$\mathcal{H}(\mathcal{C}):=\{ H(C)=[n]\setminus C:\, C\in  \mathcal{C}_L\}$$ 
be the associated family of hyperplanes of $M^\star$.
We recall  that a pair $\{X, Y\}$ of subsets of the ground set  $[n]$ is a 
{\em modular pair}
of $M([n])$ if
$$\operatorname{r}(X)+ \operatorname{r}(Y)=\operatorname{r}(X\cup Y)+
\operatorname{r}(X\cap Y).$$
\begin{proposition}\label{dual}
Let $\{C_1, C_2\}$ be a pair of circuits of $M$ and
$\{H(C_1), H(C_2)\}$ be the associated hyperplanes of $M^\star$.
The following four conditions are equivalent:
\item [$\,\,\,\circ$] $\{C_1, C_2\}$ is a modular pair of circuits of $M$,
\item [$\,\,\,\circ$] $\{H(C_1), H(C_2)\}$  is a modular pair
of hyperplanes of $M^\star$,
\item [$\,\,\,\circ$] $\operatorname{r}_M(C_1\cup C_2)=|C_1\cup C_2|-2$,
\item [$\,\,\,\circ$] $\operatorname{r}_{M^\star}(H(C_1)\cap H(C_2))
=\operatorname{r}({M^\star})-2(=n-r-2)$.\qed
\end{proposition}
\begin{definition}[\cite{T}]\label{modular}
{\em  We say that the family of circuits $\mathcal{C}'$,  
$\mathcal{C}'\subseteq\mathcal{C}(M)$,  is a
{\em linear class of circuits} 
if, given a modular pair of circuits  in $\mathcal{C}'$, all the 
circuits contained in the union of the modular pair are also in $\mathcal{C}'$.
}
\end{definition}
In the following we will always denote by $\mathcal{C}_L$  a linear 
class of circuits of the matroid $M$.
\begin{definition}\label{modular1}
{\em
We say that the family $\mathcal{H}$ of hyperplanes of $M$  is  a 
{\em  linear class of hyperplanes     of $M$} if, given a modular pair of 
hyperplanes  in $\mathcal{H}$, all the hyperplanes of $M$ containing the 
intersection of the pair are also in  $\mathcal{H}$.}
\end{definition}
%%%%%%%%%%
The following corollary is a direct consequence of
Proposition~\ref{dual} and  Definitions \ref{modular} and \ref{modular1}.
\begin{corollary}
The following two assertions are equivalent:
\item [$\,\,\,\circ$]  The family\, $\mathcal{C}'$   is a linear class of circuits of $M$;
\item [$\,\,\,\circ$]  The set $\mathcal{H}(\mathcal{C}') $  is a linear class of hyperplanes of $M^\star$.\qed
\end{corollary}
\begin{remark}
{\em

The linear class of hyperplanes $\mathcal{H}(\mathcal{C}_L)$ of $M^\star$ determines a 
single-element extension $$M^\star([n])\stackrel{\mathcal{H}(\mathcal{C}_L)}
\hookrightarrow N^\star([n+1]),$$
where  $\{n+1\}$ is in the closure in $ N^\star([n+1])$ of a hyperplane $H$ of $M^\star([n])$,
if and only if $H\in \mathcal{H}(\mathcal{C}_L)$. 
Two special cases   occur:
\begin{itemize}
\item[$\circ$] If  $\mathcal{C}_L=\mathcal{C}(M)$ the element  $n+1$ is a coloop  of 
$ N([n+1])$.
\item[$\circ$] If $\mathcal{C}_L=\emptyset=\mathcal{H}(\mathcal{C}_L)$ the element  $n+1$ is a  is in general position in   $ N^\star([n+1]$. 
\end{itemize}

In the literature 
 $N([n+1])$ is called the 
\emph{extended
lift of $M([n])$ (determined by the linear class of circuits $\mathcal{C}_L)$}. 
}
\end{remark}
%%%%%
\begin{lemma}\label{C(N)} 
Let $N=N([n+1])$ be the extended
lift of $M([n])$ determined by the linear class of circuits\, $\mathcal{C}_L, \,\mathcal{C}_L\not = \emptyset, \mathcal{C}(M)$.
Then $N$ has the family
of circuits:
$$\mathcal{C}(N)=
\begin{cases}
\mathcal{C}_L\cup\mathcal{C}_1&\text{if}\,\,\,\, |\cup_{C\in \mathcal{C}_L}C|-\operatorname{r}_M(\cup_{C\in \mathcal{C}_L}C)=n-r-1;\\
\mathcal{C}_L\cup\mathcal{C}_1\cup \mathcal{C}_2&\,\text{otherwise,}
\end{cases}
$$
where
$$\mathcal{C}_1:=
\big\{C\cup \{n+1\}: C\in \mathcal{C}(M)\setminus \mathcal{C}_L\big\},
$$
$$
\mathcal{C}_2:=
\{C'\cup C'':  C', C''\,\, \text{is  a modular pair    of}\,\, \,\mathcal{C}(M)
\setminus \mathcal{C}_L\}.
$$
\end{lemma}
\begin{proof} The matroid $N^\star([n+1])$ has the family of hyperplanes:
$$
\mathcal{H}(N^\star)=\begin{cases}
%\begin{align}
\mathcal{H}_0\cup \mathcal{H}_1 &
\text{ if}\,\, \operatorname{r}_{M^\star}(\bigcap_{C\in \mathcal{C}_L} H(C))=1;\\
\mathcal{H}_0\cup \mathcal{H}_1\cup \mathcal{H}_2 &\,\,\text{otherwise},
%\end{align}
\end{cases}
$$
where
$$
\mathcal{H}_0:=\{H\cup\{n+1\}
 :\, H \in \mathcal{H}( \mathcal{C}_L)\},
$$
$$
 \mathcal{H}_1:=\{H(C'):\, C'\in \mathcal{C}(M)\setminus \mathcal{C}_L\},
$$
$$
 \mathcal{H}_2:=\{H'\cap H''\cup\{n+1\}: H', H'' \text{ is a modular pair   of}\, \,\mathcal{H}(\mathcal{C}(M)
\setminus \mathcal{C}_L)\}.
$$
\end{proof}
%%%%%%%%%%%%%%%%%%%%%

%%%%%%%%%%%%%%%%%%%%%%%%%%%%
\section{A bias algebra}
The pair $(M, \mathcal{C}_L)$ can be seen
as a matroidal generalization of the pair $(G, \mathcal{C}_L)$ 
( defining a biased graph) where $G$ is a  graph and  
$ \mathcal{C}_L$ a set of balanced circuits of $G$.  
A biased graph is a graph together with a linear class of circuits which 
are called balanced. It is a generalisation
of signed and gain graphs which are related to some
special class of hyperplane arrangements. In the classical graphic 
hyperplane arrangements, a hyperplane has equation of the form 
$x_i=x_j$. In the ``signed graphic'' arrangements, the equations can 
be of the form $x_i=\pm x_j$. In the ``gain graphic'' arrangements, 
the equations can be of the form   $x_i=g x_j $ (in the biased case) 
or of the form $x_i= x_j +g$ (in the lift case). All these definitions 
due to T. Zaslavsky are very natural and produce a nice theory \cite{Z, Z1} 
in connection with graphs, matroids and arrangements.
The  following bias algebra is close related to the  biased graphs 
(and its matroidal generalizations).
\begin{definition}\label{algebra2}
{\em
Let $\mathcal{C}_L$ be a linear class of circuits of the matroid $M([n])$ and 
$N=N([n+1])$ be the extended
lift of $M([n])$ determined by $\mathcal{C}_L$.
Let $\OS(N)$ be the Orlik-Solomon $\mathbf{K}$-algebra of the matroid $N$.
The {\em bias  $\mathbf{K}$-algebra} of the pair $(M, \mathcal{C}_L)$, 
denoted $ \operatorname{Z}(M, \mathcal{C}_L)$, is the graded quotient 
of the Orlik-Solomon algebra $\OS(N)$ by the  two-sided ideal  
generated by $e_{n+1}$, i.e.,
$$ 
\text{Z}(M, \mathcal{C}_L):=\OS(N)/\langle e_{n+1}\rangle.
$$
}
\end{definition}
%%%%%%%%%%%%%%%%%%%%%
\begin{remark}\cite{N}
{\em
This algebra is also known as the Orlik-Solomon algebra of the pointed matroid
$N$, with basepoint $n+1$, see \cite [Definition 3.2]{Falk}.  If $N$ may be realized by a complex hyperplane
arrangement, then $Z(M,\mathcal{C}_L)$ is isomorphic to the cohomology ring of
the complement of the decone of this arrangement with respect to the
$(n+1)^\text{st}$ hyperplane, \cite[Corollary 3.57]{OT}.
        Two special cases   occur when $M$ itself is
realizable  and $\mathcal{C}_L$ is either all of $\mathcal{C}(M)$ or the empty
set.  Indeed, suppose that $M$ is the matroid associated to a complex hyperplane
arrangement $\mathcal{A}$.  Then $Z(M,\mathcal{C}(M))$ is isomorphic to the
cohomology of the complement of $\mathcal{A}$ (i.e., the Orlik-Solomon algebra of $M$), and $Z(M,\emptyset)$ is
isomorphic to the cohomology of the complement of the affine arrangement
attained by translating each of the hyperplanes of $\mathcal{A}$ some distance
away from the origin, so that every dependent set will have empty
intersection.
}
\end{remark}
%%%%%%%%%%%%%
\begin{theorem}\label{order}
The  bias $\mathbf{K}$-algebra $ \operatorname{Z}(M, \mathcal{C}_L)$ 
is independent of the order of the elements
of $M([n])$, i.e., it is an invariant of the pair $(M, \mathcal{C}_L)$.
For every linear class\, $\mathcal{C}_L$,
 the algebra $ \operatorname{Z}(M, \mathcal{C}_L)$ is isomorphic to the quotient of
the exterior
$\mathbf{K}$-algebra
\begin{equation}\label{E}
\mathcal{E}:=\bigwedge\big(\bigoplus_{i=1}^{n}\mathbf{K} e_i\big)
\end{equation}
by the  two-sided ideal $\langle \I(\mathcal{C}_L)\rangle $ generated by the set
$$\I(\mathcal{C}_L):=
\big\{\partial e_C:\, C\in \mathcal{C}_L, |C|>1\big\}\cup 
\big\{e_C: C\in \mathcal{C}(M)\big\}.
$$
\end{theorem}
\begin{proof}
Since the Orlik-Solomon  $\mathbf{K}$-algebra  $\OS(N)$ does not depend 
of the ordering of the ground set the first part of the theorem follows.
The second assertion is a straightforward consequence of 
Lemma~\ref{C(N)}. 
\end{proof}
As the element $e_{n+1}$ does not appear in the algebra 
$ \operatorname{Z}(M, \mathcal{C}_L)$ we will omit  it.
We remark that the monomial $e_X$, $X\subseteq [n]$, in 
$ \operatorname{Z}(M, \mathcal{C}_L)$  is different from zero
if and only if $X$ is an independent set of $M$.
\begin{corollary}
The bias $\mathbf{K}$-algebra $ \operatorname{Z}(M, \mathcal{C}(M))$ 
is  the  Orlik-Solomon\newline $\mathbf{K}$-algebra of \,$\OS(M)$. Furthermore the bias $\mathbf{K}$-algebra $ \operatorname{Z}(M, \emptyset)$
is isomorphic to the quotient of
the exterior algebra $(\ref{E})$ by the  two-sided ideal  generated by the set
$\{e_C: C\in \mathcal{C}(M)\}
$.\qed
\end{corollary}

\begin{corollary} If\, $\operatorname{cl}_{M'}(n+1)=n+1$, 
the bias $\mathbf{K}$-algebra $ \operatorname{Z}(M, \mathcal{C}_L(M))$ is the quotient of the exterior algebra $(\ref{E})$
 by the  two-sided ideal generated by the set
$$
\big\{\partial e_C:\, C\in \mathcal{C}_L, |C|>1\big\}\cup 
\big\{e_C: C\in \mathcal{C}(M)\big\}.
$$
\qed
\end{corollary}
%%%%%%%%%%%%%%%%%%%%%%%
\begin{definition}
{\em Given an independent set $I$,  
a non-loop element $x\in \operatorname{cl}(I)\setminus I$
is said to be
{\em $\mathcal{C}_L$-active in $I$} if $C( x,I)$ (i.e., the unique circuit contained
in $I\cup x$) is a circuit of the family $ \mathcal{C}_L$ and $x$ 
is the smallest element of $C( x,I)$.  An independent set  
with at least one $\mathcal{C}_L$-active element is said to be 
{\em $\mathcal{C}_L$-active}, and {\em $\mathcal{C}_L$-inactive} otherwise. 
We denote by $a(I)$ the smallest 
$\mathcal{C}_L$-active element in an active independent set  $I$. 
}
\end{definition}
%%%%%%%%%%%%%%%%%%%
\begin{definition}
{\em We say that a subset $U\subseteq [n]$  is a 
$\mathcal{C}_L$-\emph{unidependent} 
(set of  $M$) if it contains a unique circuit $C(U)$ of $M$, $C(U)\in\mathcal{C}_L$ and $|C(U)|>1$.\par
We say that a $\mathcal{C}_L$-unidependent  set 
$U$ is  $\mathcal{C}_L$-\emph{inactive} if the minimal element of
$C(U)$, $\operatorname{min} 
C(U)$, is the the smallest $\mathcal{C}_L$-active element
of the independent set $U\setminus\operatorname{min} C(U)$. 
Otherwise the set 
$U$ is said $\mathcal{C}_L$-\emph{active}.
}
\end{definition}
%%%%%%%%%%%%%%%
\begin{definition}
{\em  For every circuit 
$C\in \mathcal{C}_L,$ $|C|>1$,  the
  set $C\setminus \operatorname{min}(C),$\, is said to be a
$\mathcal{C}_L$-\emph{broken circuit}. 
The family  
of  $\mathcal{C}_L$-inactive independents, denoted 
$\operatorname{NBC}_{\mathcal{C}_L}$,  is the family of independent 
sets of $M$  not containing a $\mathcal{C}_L$-broken circuit.}
\end{definition}
%%%%%%%%%%%%%%%%
Set $$\mathbf{nbc}_{\mathcal{C}_L}:=\{e_I:\, 
I\in\operatorname{NBC}_{\mathcal{C}_L}\}\,\,\, \text{and}$$
\begin{equation*}\mathbf{b}_{\I(\mathcal{C}_L)}:=\{\partial e_U: 
U \text{ is $\mathcal{C}_L$-inactive unidependent\,}\}\cup 
\end{equation*}
\begin{equation*}
\cup \{e_D:  D \text{ is dependent}\}.
\end{equation*}
\begin{theorem}\label{basis}
The sets $\mathbf{nbc}_{\mathcal{C}_L}$ and  $\mathbf{b}_{\I(\mathcal{C}_L)}$\, are bases, respectively of the  bias $\mathbf{K}$-algebra $ \operatorname{Z}(M, \mathcal{C}_L)$ and  of the ideal $\langle \I(\mathcal{C}_L)\rangle$.
\end{theorem}
\begin{proof}
We will show the two statements at the same time by proving that both
sets are spanning and that they have the correct size.
Let  $I$  be an independent set of M. If $I$ is 
$\mathcal{C}_L$-active then we have
$$e_I=\sum _{x\in C(a (I) ,I)\setminus a(I)} \zeta_x\, 
e_{I\cup a(I) \setminus x},$$ where $\zeta(x)\in \{-1,1\}$. This 
is an expression for $e_I$ whit respect to lexicographically smaller 
$e_X$ where $X$ is an independent of $M$ and $|X|=|I|$. 
By induction, we get that the set $\mathbf{nbc}_{\mathcal{C}_L}$ is  a 
generator of the  graded  algebra $\operatorname{Z}(M,\mathcal{C}_L)$.\par
Let $U$ be a $\mathcal{C}_L$-unidependent set of $M$. Suppose that 
$U$ is $\mathcal{C}_L$-active  and let $a=\operatorname{min} C(U)$  and set
$I:=C(U)\setminus a$. Note that $\{C(U), C(a(I), I)\}$ is 
a modular pair   of circuits  of $\mathcal{C}_L$,
so every circuit contained in the cycle $C(U)\cup C(a(I), I)$ 
is in $\mathcal{C}_L$. From the definition of the map $\partial$
we know
that
$$\partial e_U=\sum _{x\in C(U)\setminus a} \epsilon_x\, 
\partial e_{U\cup a(I) \setminus x},$$ where 
$\epsilon_x\in \{-1,1\}$. This 
is an expression for $\partial e_U$ with respect to lexicographically smaller 
$\partial e_X$, where $X$ is a $\mathcal{C}_L$-unidependent and $|U|=|X|$.
By induction, we get that the set  $\mathbf{b}_{\I(\mathcal{C}_L)}$  
is    a generator of $\langle \I(\mathcal{C}_L)\rangle$.
By the definition of $\operatorname{Z}(M,\mathcal{C}_L)$, we know that 
$$ \operatorname{dim}(\operatorname{Z} (M,\mathcal{C}_L))+ 
\operatorname{dim}(\langle \I(\mathcal{C}_L)\rangle)= \operatorname{dim}
(\mathcal{E})=
2^n.$$ Given
a subset $X$ of $[n]$,  it is either dependent or independent 
$\mathcal{C}_L$-active or
independent $\mathcal{C}_L$-inactive. To every independent $\mathcal{C}_L$-active  independent set $I$
 corresponds uniquely the
 unidependent $\mathcal{C}_L$-inactive $I\cup a(I)$. We have then that $$|\mathbf{nbc}_{\mathcal{C}_L}(M)|+
|\mathbf{b}_{\I(\mathcal{C}_L)}|=2^n.$$
\end{proof}
%%%%%%%%%%%%%%%%%%%%%%%%%%%%%
We define the deletion  and contraction 
operation for an arbitrary subset of circuits 
$\mathcal{C}^{\prime}\subseteq \mathcal{C}(M)$ setting:
$$\mathcal{C}^{\prime}\setminus x:=\{ C\in \mathcal{C}^{\prime}
:\,  x\not \in C\}$$
and $$\mathcal{C}^{\prime}/ x:=
\begin{cases}
\mathcal{C}^{\prime}\setminus x\,\,\,\,\, \text{if $x$ is a loop of $M$,}\\
\big\{C\setminus x:\, x\in C 
\in\mathcal{C}^{\prime}\big\}\uplus\big\{C\in 
\mathcal{C}^{\prime}:\,\, x\not\in \operatorname{cl}_M(C)\big\}\,\,\,\, \,\text{otherwise}.
\end{cases}$$ 
From the preceding definition, we can see that given a circuit
$C$ of $\mathcal{C}^{\prime}/ x$, where $x$ is  a non-loop of $M$, there exists a unique circuit 
$\widehat{C}\in \mathcal{C}^{\prime}$ such that
$$\widehat{C}:=
\begin{cases}C\cup x\,\,\, \text{if}\,\,\, x\in \cl_M (C),\\
C\,\,\, \text{otherwise}.
\end{cases}
$$
%%%%%%%%%%%%%%%%%%%%%%%%%%%%%
%%%%%%%%%%%%%%%%%%%%%%%%%%%%%
\begin{proposition}
Let $M$ be a matroid and\, $\mathcal{C}_L$ be a linear class of circuits of $M$.
For an element $x$ of the matroid, the circuit sets   
$\mathcal{C}_L\setminus x$ and $\mathcal{C}_L/x$ are linear classes
of $M\setminus x$ and $M/x$, respectively.
\end{proposition}
%%%%%%%%%%%%%%%%
\begin{proof}
The statement for the deletion is clear. If $x$ is a loop the result 
is also clear for the contraction. Suppose that 
  $x$ is a non-loop of $M$.   If $Y\subseteq X$ are sets  
such that  $ \operatorname{r}_M(X)= \operatorname{r}_M(Y)+1$ then 
we have 
\begin{equation}\label{rank}
 \operatorname{r}_{M/x}(X\setminus x)= \operatorname{r}_{M/x}(Y\setminus x)
+\epsilon,\,\,\epsilon\in\{0,1\}.
\end{equation} 
So, if $\{C_1, C_2\}$ is a modular pair of circuits 
of $\mathcal{C}_L/x$, 
 $\{\widehat{C_1}, \widehat{C_2}\}$ 
is  also a modular pair of circuits of $\mathcal{C}_L$. We see also from Equation \ref{rank} that if $C\subseteq C_1\cup C_2$ is a circuit
of $M/x$  then
 $\widehat{C}\subseteq \widehat{C_1}\cup\widehat{C_2}$,
so $\widehat{C}\in \mathcal{C}_L$ and necessarily $C\in \mathcal{C}_L/x$.
\end{proof}
%%%%%%%%%%%%%%%%
\begin{definition}
For a pair $(M,\mathcal{C}_L)$ and  an element $x$  of
$M$,
we define the deletion and  the contraction of the pair
$(M,\mathcal{C}_L)$ by:
$$(M, \mathcal{C}_L)\setminus x:= (M\setminus x,\mathcal{C}_L\setminus
x)$$
and
$$(M, \mathcal{C}_L)/ x:= (M/ x,\mathcal{C}_L/ x).$$
\end{definition}
%%%%%%%%%%%%%%%%%%%%
As a corollary of Theorem~\ref{order} we have:
%%%%%%%%%%%%%%%%%%%%
\begin{proposition}\label{contraction}
For every  element $x$ of $M,$
there is a unique monomorphism of\, vector spaces,\,
$$\mathfrak{i}_x:  \operatorname{Z}(M,\mathcal{C}_L)\setminus x\to
\operatorname{Z}(M,\mathcal{C}_L),$$ such that, for every 
independent set $I$ of $M\setminus x$, we have
$\mathfrak{i}_x(e_I)=e_I.$  \qed
\end{proposition}
%%%%%%%%%%%%%%%%%%%
\begin{proposition}
For every non-loop element $x$ of $M,$
there is a unique epimorphism of\, vector spaces,\,
$\mathfrak{p}_x:   \operatorname{Z}(M,\mathcal{C}_L)\to
 \operatorname{Z}(M,\mathcal{C}_L)/x,$ such that, for every
 subset $I=\{i_1,\ldots,i_\ell\}\subseteq [n]$, 
\begin{equation}\label{p_x}
\boldsymbol{\mathfrak{p}}_x e_I:=
\begin{cases}
\vspace{2mm}
 e_{I\setminus x} & \hspace{1mm} \mbox{\em if} \hspace{2mm}  x\in I, \\
\pm e_{I\setminus y} & \hspace{2mm} \mbox{\em if} \hspace{2mm}  \exists y\in I 
 \hspace{2mm}\mbox{\em such that
 $\{x,y\}\in \mathcal{C}_L$},\\
0 &  \hspace{2mm}\mbox{\em otherwise}. 
\end{cases}
\end{equation}
More precisely the value of the coefficient $\pm 1$ in the second case is
the sign of the permutation obtained by replacing $y$ by $x$ in $I$.
\end{proposition}
%%%%%%%%%%%%%%%%%%%%
\begin{proof}
From Theorem  \ref{order},  it is enough to
prove that the map $\mathfrak{p}_x$ is well determined, i.e., for all 
$\mathcal{C}_L$-unidependent $U=(i_1,\dotsc,i_{m})$
  set of $M$, we have $$\mathfrak{p}_x \partial e_U=0\in 
\I(\mathcal{C}_L/x).$$ 
We can also suppose that 
$x$ is the last element $n$.  
Note that if $n\in U$ then 
${{U}}\setminus n$ is a $\mathcal{C}_L/n$-unidependent set of
$M/n$. If $n\not \in U$ but there is $y\in U$ and $\{n, y\}\in 
\mathcal{C}_L$, we know that 
$e_U=\pm e_{U\setminus y\cup n}$ in  $ \operatorname{Z}(M,\mathcal{C}_L)$. 
Suppose that  $n\not\in U$  and that there does not  exist $y\in U$ such that 
$\{n, y\}\in \mathcal{C}_L$. Then it is clear that
$\mathfrak{p}_n {\partial}e_{U}=0$.
Suppose that  $n\in U$. It is easy to see that
$$\pm\mathfrak{p}_n \partial e_U=\sum_{j=1}^{m-1} e_{U\setminus\{j,n\}}=0.$$
Finally, if  an independent set $I$ of $M$ contains an element $y$ 
such that  $\{x,y\}$ is a circuit 
in $\mathcal{C}_L$, 
we know that there is a scalar $\chi(I; x,y)\in \{-1,1\}$ 
such that $e_I=\chi(I; x,y)e_{I\setminus y\cup x}.$
More precisely the value of $\chi(I; x,y)\in \{-1,1\}$ is
the sign of the permutation obtained by replacing $y$ by $x$ in $I$.
\end{proof}
%%%%%%%%%%%%%%%%%%%%%%
\begin{theorem}\label{sequence}
Let $M$ be a loop free matroid  and $\mathcal{C}_L$ be a linear 
class of circuits of $M$.
For every element $x$ of $M$, there is a splitting
short exact sequence of vector spaces
\begin{equation}\label{exact}
    0\to  \operatorname{Z}(M,\mathcal{C}_L)\setminus x\stackrel{\mathfrak{i}_x}\longrightarrow 
    \operatorname{Z}(M, \mathcal{C}_L)\stackrel{\mathfrak{p}_x}\longrightarrow 
    \operatorname{Z}(M,\mathcal{C}_L)/x\to 0.
    \end{equation}
\end{theorem}
%%%%%%%%%%%%%%%%%
\begin{proof}
From the definitions we know that $\mathfrak{p}_x\circ\mathfrak{i}_x$, 
is the null map
so $\mathrm{Im}(\mathfrak{i}_x)\subseteq \mathrm{Ker}(\mathfrak{p}_x).$ 
We will prove the equality  $\dim(\mathrm{Ker}(\mathfrak{p}_n))=\dim
(\mathrm{Im}(\mathfrak{i}_n)).$ By a reordering  of the elements of $[n]$
we can suppose that  $x=n.$
The minimal  $\mathcal{C}_L/n$-broken circuits of $M$
are the minimal sets $X$ such that either $X$
or $X\cup \{n\}$ is a $\mathcal{C}_L$-broken circuit of $M$ 
(see    \cite[Proposition 3.2.e]{Bry}).  Then
%%%%%%%%%
\begin{equation*}
\operatorname{NBC}_{\mathcal{C}_L/n}=\big\{X: X\subseteq 
 [n-1]~~\mbox{and}~~X\cup\{n\}\in \operatorname{NBC}_{\mathcal{C}_L}\big\}
 \end{equation*}
%%%%%%
and we have
\begin{equation}\label{NBC}
\operatorname{NBC}_{\mathcal{C}_L}=
\operatorname{NBC}_{\mathcal{C}_L\setminus n}\uplus
\big\{I\cup
n: I\in
\operatorname{NBC}_{\mathcal{C}_L/n}\big\}.
\end{equation}
%%%%%%%
So $\dim(\mathrm{Ker}(\mathfrak{p}_n))=\dim
(\mathrm{Im}(\mathfrak{i}_n)).$ There is a  
 morphism of vector spaces 
$$\mathfrak{p}^{-1}_n:   \operatorname{Z}(M, \mathcal{C}_L)/n\to
  \operatorname{Z}(M, \mathcal{C}_L),$$ where, for every $I\in 
\operatorname{NBC}_{\mathcal{C}_L/n}$, we have
  $\mathfrak{p}^{-1}_n e_I:= e_{I\cup n}$.
It is clear that $\mathfrak{p}_n\circ \mathfrak{p}^{-1}_{n}$ is the
identity map. From Equation~(\ref{NBC}) we conclude that the exact
sequence~(\ref{exact}) splits.
\end{proof}

\begin{remark}
{\em
A large class of algebras, the so called 
$\chi$-algebras (see \cite{CF} for more details),
contain the    Orlik-Solomon, Orlik-Terao \cite{OT1} 
(associated to vectorial matroids) and Cordovil algebras \cite{Cor} 
(associated to oriented matroids).
Following the same ideas it is possible to generalize the definition
of the bias algebras  and obtain a class of bias $\chi$-algebras,  
determined by
a pair $(M,\mathcal{C}_L)$, and that contain all the mentioned algebras. 
}
\end{remark}
%%%%%%%%%%%%%%%%%%%%%%%%%%
%%%%%%%%%%%%%%%%%%%%%%%%%%
%%%%%%%%%%%%%%%%%%%%%%%%%%

%%%%%%%%%
Similarly to  \cite{CF}, we now construct, making use  
of iterated contractions,  the dual basis 
$\boldsymbol{nbc}_{\mathcal{C}_L}^*$ of the standard basis  
$\boldsymbol{nbc}_{\mathcal{C}_L}$.
Let $\operatorname{Z}(M, \mathcal{C}_L)_h$ be the subspace of
$\operatorname{Z}(M, \mathcal{C}_L)$ generated by the set 
$$\{e_X: X\,\,\, \text{is an independent set of}\,\,\, M\,\,\, 
\text{and}\,\,\, |X|=h\}.$$
We associate to  the (linearly ordered)  independent set 
$I=(i_{1},\dotsc,  i_h)$ of
$M$ the linear form on 
$\operatorname{Z}(M, \mathcal{C}_L)_h,$\, 
$\boldsymbol{\mathfrak{p}}_{I}:
\operatorname{Z}(M, \mathcal{C}_L)_h\to  \K,$
%%%%%%%
\begin{equation}\label{resid}
\boldsymbol{\mathfrak{p}}_{I}:=
\boldsymbol{\mathfrak{p}}_{e_{i_{
1}}}\!\circ\boldsymbol{\mathfrak{p}}_{e_{i_{2}}}\!\circ
\cdots\,\circ
\boldsymbol{\mathfrak{p}}_{e_{i_h}}.
\end{equation}
%%%%%%%
 We also associate to the linearly ordered independent 
$I=(i_{1},\dotsc,  i_j)$ 
the flag of its final independent subsets, defined by $$\{I_t: I_t=
(i_{t}, \ldots,i_j),\, 1\le t\le j\}.$$
%%%%%%%%%%%%%%%%%%%%%%
\begin{proposition}\label{resi}
Let $I=(i_{1},\dotsc,  i_h)$ and $J=(j_1,\ldots,j_h)$
be two linearly ordered independents of $M$, then we have 
$\boldsymbol{\mathfrak{p}}_{I}(e_{J})\not=0$ if and only if there is a
permutation $\tau \in \os_h$  such that for every $1\le t\le h$,
$j_{\tau(t)}\in\cl (I_t)$ and $C(j_{\tau(t)},I_t)\in \mathcal{C}_L$. 
When the permutation $\tau$ exists, it is unique and we have 
$\boldsymbol{\mathfrak{p}}_{I}(e_{J})=\operatorname{sgn}(\tau)$.
In particular we have
$\boldsymbol{\mathfrak{p}}_{I}(e_{I})=1$ for any independent set $I$.
\end{proposition}
%%%%%%%%%%%%%%%%%
\begin{proof}
The first equivalence is very easy to prove in both directions.
To obtain the expression of $\boldsymbol{\mathfrak{p}}_{I}(e_{J})$
we just need to iterate $h$ times the formula of contraction of Proposition 
\ref{contraction}.
With the definition of the permutation $\tau$ we know that $
\boldsymbol{\mathfrak{p}}_{I}(e_{\tau(1)}\wedge\cdots\wedge e_{\tau(h)})=1$.
By the antisymmetric of the wedge product we also have that
$e_J=\sgn (\tau)\times e_{\tau(1)}\wedge\cdots\wedge e_{\tau(h)}$.
And finally the
last result comes from the fact that if $I=J$ then clearly $\tau=\id.$
\end{proof}
%%%%%%%%%%%%%%%%%%%%%%%%%%
\begin{theorem}\label{dbasis} 
The set 
\,$\{\boldsymbol{\mathfrak{p}}_{I}: 
I\in \operatorname{NBC}_{\mathcal{C}_L}\}$ is the dual 
basis  of the standard basis\, $\mathbf{nbc}_{\mathcal{C}_L}$ of\,  $ \operatorname{Z}(M, \mathcal{C}_L)$. 
\end{theorem}
%%%%%%%%%%%
\begin{proof} 
Pick two elements $e_I$ and $e_J$ in
$\mathbf{nbc}_{\mathcal{C}_L}$, $|I|=|J|=h$. We just need to prove that
$\mathfrak{p}_{I}(e_{J})=\delta_{IJ}$ (the Kronecker delta). From 
the preceding proposition we already have that $\mathfrak{p}_{I}(e_{I})=1$.
Suppose for a contradiction that there exists a permutation
$\tau$ such that 
$j_{\tau(t)}\in\cl (I_t)$ and $C(j_{\tau(t)},I_t)\in \mathcal{C}_L$
for every $1\le t\le h$.
Suppose that $j_{\tau(m+1)}=i_{m+1}, 
\dotsc, j_{\tau(h)}=i_h$ and $i_m\not =j_{\tau({m})}.$ 
Then there is a circuit $C\in {\mathcal{C}_L}$ 
such that $$i_m, j_{\tau(m)}\in 
C\subseteq \{i_m, j_{\tau(m)}, i_{m+1}, i_{m+2}, \dotsc
,i_h\}.$$
If $j_{\tau(m)}<i_m$ [resp. $i_m<j_{\tau(m)}$] we conclude that $I\not \in 
\text{NBC}_{\mathcal{C}_L}$ 
 [resp. $J\not \in \text{NBC}_{\mathcal{C}_L}$],  a 
contradiction.
\end{proof}
%%%%%%%
The following corollary is an extension of results 
of \cite{cor}, \cite{Cor} and \cite{CF}.
%%%%%%%%
\begin{corollary}\label{that}
Let $J=\{j_1,\ldots,j_\ell\}$ be an independent set of 
$M$ such that the expansion of
$e_J$ in $\mathbf{nbc}_{\mathcal{C}_L}$ is\, $e_J=
\sum_{I\in \mathbf{nbc}_{\mathcal{C}_L}}
\xi (I,J) e_I$.
%%%%%%%%%% 
Then the following  are equivalent:
\item [$\,\,\,\circ$] $\xi (I,J)\not =0,$
\item [$\,\,\,\circ$] there exists a permutation $\tau$ such that
$e_{\tau(t)}\in\cl(I_t)$ and $C(j_{\tau(t)},I_t)\in \mathcal{C}_L$
for every $1\le t\le h$.
%%%%%%%%%%%
Moreover, in the case where  $\xi (I,J)\not =0$ we have 
$\xi(I,J)=\sgn(\tau).$\qed
\end{corollary}
%%%%%%%%%% 


\begin{thebibliography}{999}
\bibitem{Bry} Brylawski, T.:~The broken-circuit complex.
{\em Trans. Amer. Math. Soc.} \textbf{234} (1977), no. 2, 417\---433.
\bibitem{cor}Cordovil, R., and  Etienne, G.:~A note on the Orlik-Solomon 
algebra. \emph{European J. of Combin.} \textbf{22} (2001), 
165\---170.

\bibitem{Cor} Cordovil, R.:\,\,A commutative algebra for oriented matroids.
{\em Discrete and Computational Geometry} \textbf{27} (2002), 73\---84. 

\bibitem{CF} Cordovil, R. and Forge, D.:\,\, Diagonal bases in Orlik-Solomon 
type algebras. {\em Annals of Combinatorics} \textbf{7} (2003), 247\---257.
\bibitem{Falk} Falk, Michael J.:\,\,Combinatorial and algebraic structure 
in Orlik\---Solomon algebras.
{\em European J. Combin.} \textbf{22}, (2001), no. 5, 687\---698.


\bibitem{OS1} Orlik, Peter; Solomon, Louis:~Combinatorics and topology 
of complements of hyperplanes. \emph{Invent. Math.} 
\textbf{56} (1980), no. 2, 167\---189. 

\bibitem{OT} Orlik, Peter; Terao, Hiroaki:\,\,Arrangements of
Hyperplanes. Grun\-dlehren der Mathematischen
Wissenschaften [Fundamental Principles of Mathematical Sciences],
300.  {\em Springer-Verlag, Berlin}, 1992.

\bibitem{OT1}Orlik, Peter; Terao, Hiroaki:~Commutative algebras 
for arrangements.
 \emph{Nagoya Math. J.}
\textbf{134} (1994), 65\---73. 


\bibitem{OT2} Orlik, Peter; Terao, Hiroaki:\,\,Arrangements and 
hypergeometric integrals.  MSJ Memoirs \textbf{9}. \emph{Mathematical 
Society of
Japan}, Tokyo, 2001.

\bibitem{T} Tutte, W. T.:\,\,  Lectures on matroids. \emph{ J. Res. Nat. Bur. Standards Sect.} B  \textbf{69B}  (1965), 1\---47.
\bibitem{N} Proudfoot, Nicholas, private communication.
 \bibitem{Z} Zaslavsky, T.:\,\, Biased graphs. I.  Bias, balance, and gains. {\em J. Combin. Theory Ser. B} \textbf{47} (1989), 32\---52.

 \bibitem{Z1} Zaslavsky, T.:\,\, Biased graphs.  II.  The three matroids.
{\em J. Combin. Theory Ser. B} \textbf{51} (1991), 46\---72.


\end{thebibliography}
\end{document}